\documentclass{amsart}

\newtheorem{theorem}{Theorem}[section]

\theoremstyle{definition}

\theoremstyle{remark}

\begin{document}

\title[\empty]{An elementary trigonometric equation}


\author{Victor H. Moll}
\address{Department of Mathematics,
Tulane University, New Orleans, LA 70118}
\email{vhm@math.tulane.edu}


\date{\today}

\keywords{Quadratic radicals, Gauss sums}

\begin{abstract}
We provide a systematic study of the trigonometric equation
\begin{equation}
A \tan a + B \sin b = C.
\nonumber 
\end{equation}
\noindent
In particular, the classical example 
\begin{equation}
\tan \frac{3 \pi}{11} + 4 \sin \frac{2 \pi}{11} = \sqrt{11},
\nonumber
\end{equation}
\noindent
appears in a natural form $\sqrt{11}$. A second proof 
involving Gaussian sums is also discussed. 
\end{abstract}

\maketitle

\newcommand{\ba}{\begin{eqnarray}}
\newcommand{\ea}{\end{eqnarray}}
\newcommand{\ift}{\int_{0}^{\infty}}
\newcommand{\nn}{\nonumber}
\newcommand{\no}{\noindent}
\newcommand{\lfour}{\ln \left(1 + 4 \theta^{4}/k^{4} \right)}

\newtheorem{Definition}{\bf Definition}[section]
\newtheorem{Thm}[Definition]{\bf Theorem} 
\newtheorem{Example}[Definition]{\bf Example} 
\newtheorem{Lem}[Definition]{\bf Lemma} 
\newtheorem{Cor}[Definition]{\bf Corollary} 
\newtheorem{Prop}[Definition]{\bf Proposition} 
\numberwithin{equation}{section}

\section{Introduction} \label{intro} 
\setcounter{equation}{0}

In the Problems and Solutions section of the September 2006 issue of 
The College Mathematics Journal \cite{rabinowitz1}, one is told of the identity
\begin{equation}
\tan \frac{3 \pi}{11} + 4 \sin \frac{2 \pi}{11} = \sqrt{11},
\label{one}
\end{equation}
\noindent
that appeared as Problem 218 of the same journal \cite{unknown}. The problem
asks to prove the related identites 
\begin{equation}
\tan \frac{\pi}{11} + 4 \sin \frac{3 \pi}{11} = \sqrt{11},
\label{two}
\end{equation}
\noindent
and 
\begin{equation}
\tan \frac{4\pi}{11} + 4 \sin \frac{\pi}{11} = \sqrt{11}.
\label{three}
\end{equation}

These evaluations have appeared in classical texts. For instance, (\ref{one}) 
is an exercise in Hobson's Trigonometry 
text \cite{hobson1} (pages 123 and 382). It also appears as Exercise $14$ on
page $270$ of T. J. I. Bromwich's treatise on Infinite Series \cite{bromwich}. 

During the period 1991-2003 it was my privilege to work with George Boros. 
First as an advisor and then as colleague. Our method of work  was rather
unorthodox. George has a very good knowledge of numbers, so he will bring 
me many pages of formulas and my role was to try to figure out where these
things fit. Around 1997 he showed me the `$\sqrt{11}$ problem'. This note 
is a reflection of what I learned from him.

\section{The reduction}
\setcounter{equation}{0}

We consider the equation
\begin{equation}
A \tan a + B \sin b = C,
\label{eqn1}
\end{equation}
\noindent
where $A, \, B$ and $C^{2}$ are {\em rational} numbers. We may assume that 
$A=1$ and we write (\ref{eqn1}) as 
\begin{equation}
\sin a + B \cos a \sin b = C \cos a. \label{eqn2}
\end{equation}
\noindent
The procedure described here will give {\em some} solutions of (\ref{eqn2}). 
In particular, the radical $\sqrt{11}$ will appear in a natural form. \\

The identity
\begin{equation}
\cos a \sin b = \tfrac{1}{2} \left( \sin(a+b) - \sin(a-b) \right),
\end{equation}
\noindent
converts (\ref{eqn2}) into
\begin{equation}
\sin a + \tfrac{1}{2}B \sin(a+b) - \tfrac{1}{2}B \sin(a-b) = C \cos a. 
\label{eqn3}
\end{equation}

\medskip

\begin{Lem}
\label{reduced}
Assume the angles $a, \, b, \, c$ satisfy (\ref{eqn3}). Then
\begin{eqnarray}
\tfrac{1}{4}(2-2C^{2}+B^{2}) & = & 
\tfrac{1}{4}(2+2C^{2}-B^{2}) \cos (2a) + \tfrac{1}{8}B^{2} \cos(2a+2b) + 
\nonumber  \\
& + & \tfrac{1}{8}B^{2} \cos(2a-2b) 
+ \tfrac{1}{2}B \cos(2a+b) + 
\tfrac{1}{4}B^{2} \cos(2b) \nonumber \\
 & - & \tfrac{1}{2}B \cos(2a-b).
\nonumber
\end{eqnarray}
\end{Lem}
\begin{proof}
Square (\ref{eqn3}) and use 
\begin{eqnarray}
\sin u \sin v & =&  \tfrac{1}{2} \cos(u-v) - 
\tfrac{1}{2} \cos(u+v), \nonumber \\
\cos u \cos v & =&  \tfrac{1}{2} \cos(u-v) + \tfrac{1}{2} \cos(u+v), \nonumber
\end{eqnarray}
\noindent
the standard formulas for double angle 
\begin{equation}
\sin^{2}u = \tfrac{1}{2}(1- \cos(2a)) \text{ and } 
\cos^{2}u = \tfrac{1}{2}(1+ \cos(2a)), 
\end{equation}
and 
\begin{equation}
\sin^{2}(a+b) + \sin^{2}(a-b) = 1 - \cos(2a) \cos(2b)
\end{equation}
to produce the result.
\end{proof}

\section{Some special cases}
\setcounter{equation}{0}

In order to produce some special solutions of the equation in 
Lemma \ref{reduced}, we begin by
matching the different coefficients appearing in it. For example, 
forcing 
\begin{equation}
\tfrac{1}{8}B^{2} = \tfrac{1}{2}B,
\end{equation}
\noindent
that is, choosing $B=4$, converts it into
\begin{eqnarray}
\tfrac{1}{2}(9-C^{2}) & = & 
\tfrac{1}{2}(C^{2}-7) \cos(2a) + 2 \cos(2a+2b) + 2 \cos(2a-2b) + 
\label{sym1} \\
&  & + 2 \cos(2a+b) + 4 \cos(2b) - 2 \cos(2a-b). \nonumber 
\end{eqnarray}
\noindent
We now simplify this equation further by imposing  the condition
\begin{equation}
\tfrac{1}{2}(C^{2}-7) = 2,
\end{equation}
\noindent
to make almost every 
coefficient on the right hand side of (\ref{sym1}) equal to 
$\pm 2$. This produces $C = \pm \sqrt{11}$. \\

The equation (\ref{sym1}) becomes
\begin{eqnarray}
-\tfrac{1}{2} & = & \cos(2a-2b) - \cos(2a-b) + \cos(2a) + \cos(2a+b) + 
 \label{special1} \\
 & + & \cos(2a+2b) + 2 \cos(2b). \nonumber
\end{eqnarray}

Observe that five out the six angles appearing in (\ref{special1}) are in 
arithmetic progression: $2a-2b, \, 2a-b, \, 2a, \, 2a+b$ and $2a+2b$. 
We now proceed to establish further restrictions on the angles 
$a$ and $b$ so that the final angle, namely $2b$, is also part of this 
progression. This will determine {\em some} solutions of the original 
problem. 

\subsection{The first example}
We assume first $2b=2a-b$, that 
is, $2a=3b$. Then (\ref{special1}) becomes
\begin{equation}
\cos(2a-2b) + \cos(2a-b) + \cos(2a) + \cos(2a+b) + \cos(2a+2b) = 
- \tfrac{1}{2}. \label{case1}
\end{equation}

This identity is simplified using the closed form expression for sums of 
cosines in arithmetic progression:
\begin{equation}
\sum_{k=0}^{n-1} \cos(x+ky) = \frac{\cos(x + (n-1)y/2) \, \sin(ny/2)}
{\sin(y/2)}. 
\label{sum-cos}
\end{equation}
\noindent
This formula is easy to establish and it can be found as $1.341.3$ in the 
table \cite{gr}. In our first case, we have 
$x=2a-2b, y=b$ and $n=5$ so that (\ref{sum-cos}) produces 
\begin{equation}
\frac{\sin(5b/2) \, \cos(2a)}{\sin(b/2)} = - \tfrac{1}{2}.
\label{sum1}
\end{equation}
\noindent
We exclude the case where $b$ is an integer multiple of $2 \pi$, since in that 
situation equation (\ref{eqn1}) is not interesting. 
Using $2a=3b$ and the elementary identity
\begin{equation}
\sin x \cos y = \tfrac{1}{2} \left( \sin(x+y) + \sin(x-y) \right),
\end{equation}
\noindent
(\ref{sum1}) produces 
\begin{equation}
\sin \left( \frac{11b}{2} \right) = 0 \text{ and } 
\sin \left( \frac{b}{2} \right) \neq 0. 
\end{equation}
\noindent
We conclude that 
\begin{equation}
b = \frac{2k \pi}{11}, \, k \in \mathbb{Z}, \, k \not 
\equiv 0 \bmod 11, 
\end{equation}
\noindent
and the angle $a$ is given by
\begin{equation}
a = \frac{3k \pi}{11}, \, k \in \mathbb{Z}, \, k \not 
\equiv 0 \bmod 11. 
\end{equation}

\medskip

Therefore, for $k \in \mathbb{Z}$ not divisible by $11$, we have found 
some solutions to (\ref{eqn1}): 
\begin{equation}
\tan \left( \frac{3k \pi}{11} \right) + 4 \sin \left( \frac{2k \pi}{11} \right) 
= \pm \sqrt{11}. 
\end{equation}
\noindent
A numerical computation of the left hand side shows that, modulo $11$, 
$k = 1, \, 3, \, 4, \, 5, \, 9$ correspond to the positive sign and 
$k = 2, \, 6, \, 7, \, 8, \, 10$ to the minus sign.  Reducing the angle to
the smallest integer multiple of $\pi/11$ we obtain the five relations 
described in the next theorem, that include (\ref{one}), (\ref{two}), and 
(\ref{three}). 

\begin{theorem}
The following identities hold
\begin{eqnarray}
\tan \left( \frac{\pi}{11} \right) + 4 \sin \left( \frac{3 \pi}{11} 
\right) & = & \sqrt{11}, \label{identity1}  \\
\tan \left( \frac{2 \pi}{11} \right) - 4 \sin \left( \frac{5 \pi}{11} 
\right) & = & -\sqrt{11}, \nonumber  \\
\tan \left( \frac{3 \pi}{11} \right) + 4 \sin \left( \frac{2 \pi}{11} 
\right) & = & \sqrt{11}, \nonumber \\
\tan \left( \frac{4 \pi}{11} \right) + 4 \sin \left( \frac{\pi}{11} 
\right) & = & \sqrt{11}. \nonumber  \\
\tan \left( \frac{5 \pi}{11} \right) - 4 \sin \left( \frac{4 \pi}{11} 
\right) & = & \sqrt{11}. \nonumber 
\end{eqnarray}
\end{theorem}

\subsection{A second example.} We now assume $a=b$. Then 
(\ref{special1}) becomes
\begin{equation}
\tfrac{1}{4}(5-C^{2}) = -\cos a + \tfrac{1}{4}(C^{2}+1) + \cos(2a)+ 
\cos(3a) + \cos(4a). 
\end{equation}
\noindent
Choose $C^{2}=3$ to make all the coefficients on the right hand side 
$\pm 1$. This yields
\begin{equation}
- \cos a + \cos(2a) + \cos(3a) + \cos(4a) = \tfrac{1}{2},
\end{equation}
\noindent 
that we write as 
\begin{equation}
1 + \cos a + \cos(2a) + \cos(3a) + \cos(4a) = \tfrac{1}{2}(3 + 4 \cos a). 
\end{equation}
\noindent
Using (\ref{sum-cos}) we obtain 
\begin{equation}
\frac{\sin(5a/2) \, \cos(2a)}{\sin(a/2)} = \tfrac{1}{2}(3 + 4 \cos a), 
\end{equation}
\noindent
that is equivalent to 
\begin{equation}
\sin(9a/2) = 2 \sin(a/2) \left[ 3 - 4 \sin^{2}(a/2) \right]. 
\label{equation2}
\end{equation}
\noindent
Now use $\sin(3t) = 3 \sin t - 4 \sin^{3}t$ to write (\ref{equation2}) as
\begin{equation}
\sin(9a/2) = 2 \sin(3a/2), \quad \text{ and } \sin(a/2) \neq 0. 
\end{equation}
Let $ x = 3a/2$ to write 
\begin{equation}
\sin x \, ( 4 \sin^{2}x - 1 ) = 0, 
\end{equation}
\noindent
that has solutions $x = \frac{\pi}{6} \times \{ 1 ,\, 5, \, 7, \, 11 \} + 
2 m \pi$. Therefore
\begin{equation}
a = \frac{\pi}{9} \times \{ 1, \, 5, \, 6, \, 7, \, 11, \, 12 \} 
+ \frac{12 m \pi}{9}. 
\end{equation}
\noindent
This yields a second 
family of solutions to (\ref{eqn1}). 

\begin{theorem}
The following identities hold
\begin{eqnarray}
\tan \left( \frac{\pi}{9} \right) + 4 \sin \left( \frac{\pi}{9} 
\right) & = & \sqrt{3}, \label{identity2} \\
\tan \left( \frac{2 \pi}{9} \right) - 4 \sin \left( \frac{2 \pi}{9} 
\right) & = & - \sqrt{3}, \nn \\
\tan \left( \frac{4 \pi}{9} \right) - 4 \sin \left( \frac{4 \pi}{9} 
\right) & = & \sqrt{3}, \nn \\
\tan \left( \frac{6 \pi}{9} \right) + 4 \sin \left( \frac{6 \pi}{9} 
\right) & = & \sqrt{3}. \nn 
\end{eqnarray}
\end{theorem}

Using the same techniques the reader is invited to check the following 
result. 

\begin{theorem}
The following identities hold
\begin{eqnarray}
\tan \left( \frac{\pi}{7} \right) - 4 \sin \left( \frac{2 \pi}{7} 
\right) & = & - \sqrt{7}, \label{identity3} \\
\tan \left( \frac{2 \pi}{7} \right) - 4 \sin \left( \frac{3 \pi}{7} 
\right) & = & - \sqrt{7}, \nn \\
\tan \left( \frac{3 \pi}{7} \right) - 4 \sin \left( \frac{\pi}{7} 
\right) & = & \sqrt{7}. \nn 
\end{eqnarray}
\end{theorem}

\section{A proof using Gaussian sums} \label{gaussian} 
\setcounter{equation}{0}

In  this section we present a proof of identity (\ref{one}) using the 
value of the Gaussian sum 
\begin{equation}
G_{n} = \sum_{j=0}^{n-1} e^{2 \pi i j^{2}/n}. 
\end{equation}
\noindent
Gauss proved that
\begin{equation}
G_{n} = \begin{cases} 
          (1+ i) \sqrt{n} & \quad \text{ if } n \equiv 0 \bmod 4, \\ 
          \sqrt{n} & \quad \text{ if } n \equiv 1 \bmod 4, \\ 
          0  & \quad  \text{ if } n \equiv 2 \bmod 4, \\ 
          i \sqrt{n} & \quad \text{ if }  n \equiv 3 \bmod 4. 
      \end{cases}
\label{gauss1}
\end{equation}
\noindent
In particular, for  $n = 11$, we have
\begin{equation}
G_{11} = i \sqrt{11}.
\end{equation}
\noindent
The reader will find a proof of (\ref{gauss1}) in \cite{apostol1} and much more 
information about these sums in \cite{berndtsums}. \\

Let $x = e^{2 \pi i/11}$, so that $x^{11}=1$. This yields
\begin{equation}
1 +x + x^{2} + \cdots + x^{9} + x^{10} = 0. 
\end{equation}
\noindent
Gauss' result gives
\begin{equation}
1 + 2( x + x^{3} + x^{4} + x^{5} + x^{9}) = i \sqrt{11}, 
\end{equation}
\noindent
where the numbers $1, \, 3, \, 4, \, 5, \, 9$ are the squares modulo $11$,
also called {\em quadratic residues}. 

Using $\sin t = (e^{it} - e^{-it})/2i$ and $\cos t = (e^{it}+e^{-it})/2$ we 
obtain
\begin{equation}
4i \sin \frac{2 \pi}{11} = 2( x - x^{-1}) = 2(x-x^{10}), 
\end{equation}
\noindent
and
\begin{eqnarray}
i \tan \frac{3 \pi}{11} & = & \frac{x^{3/2}-x^{-3/2}}{x^{3/2}+x^{-3/2}}  \
\label{identity6} \\
& = & \frac{x^{3}-1}{x^{3}+1} \nonumber \\
& = & \frac{x^{3}-x^{33}}{x^{3}+1} \nonumber \\
& = & x^{3}(1-x^{15}) \frac{1+x^{15}}{1+x^{3}} \nonumber \\
& = & x^{3}(1-x^{4}) \left[ 1 - x^{3} + x^{6} - x^{9} + x^{12} \right] 
\nonumber \\
& = & x^{3}(1-x^{4}) \left[ 1 + x - x^{3} + x^{6} - x^{9} \right] \nonumber \\
 & = & -x - x^{2} + x^{3} + x^{4} + x^{5} - x^{6} - x^{7} - x^{8} + x^{9} 
+ x^{10}. 
\nonumber 
\end{eqnarray}
\noindent
Therefore 
\begin{eqnarray}
i \tan \frac{3 \pi}{11} + 4i \sin \frac{2 \pi}{11} & = & 
x + x^{3} + x^{4} + x^{5} + x^{9} -( x^{2} + x^{6} + x^{7} + x^{8} + x^{10})
\nonumber \\
& = & 1 + 2(x + x^{3} + x^{4} + x^{5} + x^{9}) - (1 + x + x^{2} + 
\cdots + x^{10}) \nonumber \\
& = & 1 + 2 (x + x^{3} + x^{4} + x^{5} + x^{9}) \nonumber \\
& = & i \sqrt{11}. \nonumber 
\end{eqnarray}
\noindent
The proof of (\ref{one}) is complete. 

\medskip

This technique can be used to establish (\ref{two}), (\ref{three}) and 
many other identities. The reader is 
invited to check 
\begin{equation}
\tan \frac{2 \pi}{7} + 4 \sin \frac{2 \pi}{7} - 4 \sin \frac{\pi}{7} = 
\sqrt{7}, 
\end{equation}
\noindent
and 
\begin{equation}
\tan \frac{4 \pi}{19} + 4 \sin \frac{5 \pi}{19} - 4 \sin \frac{6 \pi}{19} 
+ 4 \sin \frac{9 \pi}{19} = \sqrt{19}, 
\end{equation}
\noindent
and also
\begin{equation}
\tan \frac{\pi}{9} + 2 \sin \frac{\pi}{9} - 2 \sin \frac{2 \pi}{9} + 
2 \sin \frac{4 \pi}{9} = \sqrt{3},
\end{equation}
\noindent
by these methods. 

\bigskip

\no
{\bf Acknowledgments}. The author wishes to thank T. Amdeberham for a
careful reading of the manuscript. The 
partial support of NSF-DMS 0409968 is acknowledged.

\end{document}